\documentclass[a4paper]{amsart}

\usepackage{enumerate}
\usepackage{stmaryrd}
\usepackage[all]{xy}
\usepackage{amssymb}

\newtheorem{theorem}{Theorem}[section]

\newtheorem{prop}[theorem]{Proposition}
\newtheorem{cor}[theorem]{Corollary}

\theoremstyle{definition}
\newtheorem{definition}[theorem]{Definition}
\newtheorem{example}[theorem]{Example}

\theoremstyle{remark}
\newtheorem{remark}[theorem]{Remark}

\numberwithin{equation}{section}

\newcommand{\Rr}{\mathbb R}
\newcommand{\Ss}{\mathbb S}

\renewcommand{\d}{\mathrm d}

\newcommand{\set}[1]{\left\{#1\right\}}

\newcommand{\eps}{\varepsilon}

\newcommand{\tto}{\rightrightarrows}

\newcommand{\X}{\ensuremath{\mathfrak{X}}}
\newcommand{\F}{\ensuremath{\mathcal{F}}}

\newcommand{\K}{\mathcal{K}}

\newcommand{\G}{\mathcal{G}}            % Lie groupoid
\renewcommand{\F}{\mathcal{F}}          % Fiber Lie groupoid
\newcommand{\Diff}{\text{\rm Diff}\,}   % Diffeomorphisms
\newcommand{\Ham}{\text{\rm Ham}\,}     % Hamiltonian Diffeomorphisms
\renewcommand{\Vert}{\text{\rm Vert}\,}   % Vertical subbundle
\newcommand{\Hor}{\text{\rm Hor}\,}     % Horizontal subbundle
\renewcommand{\graph}{\text{\rm graph}\,}
\newcommand{\al}{\alpha}                % section of Lie algebroid
\newcommand{\be}{\beta}                 % section of Lie algebroid
\newcommand{\Lie}{\mathcal{L}}          % Lie derivative
\renewcommand{\gg}{\mathfrak{g}}        % Lie algebra
\newcommand{\Ker}{\text{\rm Ker}\,}     % Kernel
\newcommand{\Ad}{\text{\rm Ad}\,}       % Adjoint
\newcommand{\su}{\mathfrak{su}}         % Lie algebra su(2)
\newcommand{\SU}{\text{\rm SU}}         % Lie group SU(2)

\begin{document}

\title{Poisson Fibrations and Fibered Symplectic Groupoids}

%    Information for authors
\author{Olivier Brahic}
\author{Rui Loja Fernandes}
\address{Depart.~de Matem\'{a}tica, 
Instituto Superior T\'{e}cnico, 1049-001 Lisboa, PORTUGAL} 
\email{brahic@math.ist.utl.pt, rfern@math.ist.utl.pt}
\thanks{Supported in part by FCT/POCTI/FEDER and by grants
  POCI/MAT/57888/2004 and POCI/MAT/55958/2004.}

%    General info

\date{December 15, 2006}


\begin{abstract}
We show that Poisson fibrations integrate to a special
kind of symplectic fibrations, called fibered symplectic
groupoids.
\end{abstract}
\maketitle

%%%%%%%%%%%%%%%%%%%%%%%%%%%%%%%%%%%%
%%%%%%%%%%%%%%%%%%%%%%%%%%%%%%%%%%%%
%%%%%%%%%%%%%%%%%%%%%%%%%%%%%%%%%%%%
\section{Introduction}             %
\label{sec:introduction}           %
%%%%%%%%%%%%%%%%%%%%%%%%%%%%%%%%%%%%
%%%%%%%%%%%%%%%%%%%%%%%%%%%%%%%%%%%%
%%%%%%%%%%%%%%%%%%%%%%%%%%%%%%%%%%%%

Our purpose in this paper is to explain that there is a geometric
theory of Poisson fibrations that is analogous to the theory of
symplectic fibrations. This paper makes no special claims to
originality, and improves previous works (see, e.g., Theorem \ref{main:thm})
of Vorobjev (Poisson case, \cite{Vor}) and Vaisman (Dirac case, \cite{Vai}). Our main
contribution is two folded. On the one hand, we propose an approach based
on gauge theory and Dirac geometry, which gives some natural explanations for
some of the mysterious formulas that appear in those works. On
the other hand, we look for the first time into the integration of
such structures, recovering symplectic fibrations from Poisson
fibrations. 

A symplectic fibration is a locally trivial fibration with fiber type
a symplectic manifold, admitting a collection of trivializations whose
transition functions are symplectomorphisms. Symplectic
fibrations have a long history going back to the early works of
Weinstein \emph{et al.}~\cite{Wein2,Got} and Guillemin \emph{et
al.}~(\cite{GS,GLS}). A closed 2-form on the total space of a fibration 
which restricts to a symplectic form on the fibers, determines a symplectic
fibration. Conversely, given a symplectic fibration, one may ask if there
exists a closed 2-form which is compatible with the fibration. It is 
well-known that there are non-trivial obstructions for the existence of 
such coupling forms and that one can classify all such coupling
forms. A very nice exposition of the the theory of symplectic
fibrations, where these results are discussed in detail, is given
in Chapter 6 of the monograph by McDuff and Salamon \cite{MaSa}.

We are interested in the more general notion of a \emph{Poisson fibration}, 
i.e., a locally trivial fibration with fiber
type a Poisson manifold, admitting a collection of trivializations
whose transition functions are Poisson diffeomorphisms. At first sight, the
analogous questions for Poisson geometry are either hopeless or
trivial. On the one hand, there are simple examples of fibrations
with a Poisson structure on the total space that restricts to each
fiber and which is not a Poisson fibration. On the other hand, every
Poisson fibration always admits trivially a compatible Poisson
structure: one can just declare the fibers to be Poisson
submanifolds. 

A closer inspection, however, reveals a very different point of
view. Note that for a symplectic fibration one looks for a
\emph{presymplectic structure} which intersects each fiber in a
symplectic submanifold. Therefore, for a Poisson fibration we should
look for a \emph{Dirac structure} for which each presymplectic leaf
intersects every fiber in a symplectic leaf of the Poisson structure on
the fiber (recall that a Poisson manifold is a (singular) foliation
by symplectic manifolds, while a Dirac manifold is a (singular)
foliation by presymplectic manifolds). This yields the notion of a
\emph{Dirac coupling} for a Poisson fibration. The theory of Poisson
fibrations can then be seen as a foliated version of the theory of
symplectic fibrations (this is the point of view advocated by
Vaisman \cite{Vai}).

Therefore, the questions (and answers) in the theory of Poisson
fibrations, are analogous to (and generalize) the theory of symplectic
fibrations. For example, given a Poisson fibration one would like to
know (i) if it admits coupling Dirac structures and (ii) classify all
such couplings. In this direction we have the following generalization
of a well-known result in symplectic fibrations (see 
\cite[Theorem 6.13]{MaSa}):

\begin{theorem}
\label{main:thm}
Let $p:M\to B$ be a Poisson fibration. Then the following statements
are equivalent:
\begin{enumerate}[(i)]
\item $p:M\to B$ admits a coupling a Dirac structure.
\item There exists a Poisson connection on $p:M\to B$ whose holonomy
  groups act on the fibers in a hamiltonian fashion.
\end{enumerate}
\end{theorem}

Its is also possible to classify all such coupling forms. The proof of this
result follows the same pattern as in the symplectic case: one builds
a Poisson gauge theory where coupling forms are obtained on associated
fiber bundles $M=P\times_G F$, starting from a connection on a principal
$G$-bundle and a Hamiltonian $G$-action on a Poisson manifold
$(F,\pi)$. Note that $G$ is not necessarily a finite dimensional Lie
group. 

Our second main purpose is the \emph{integration} of Poisson fibrations. 
Recall that Poisson structures are infinitesimal objects which integrate to global
objects called symplectic groupoids. There are obstructions to 
integrability which were recently understood \cite{CrFe1,CrFe2}, but we will
ignore them for the time being. Now, we will see that: 
\begin{itemize}
\item The global object associated to a Poisson fibration is a \emph{fibered
  symplectic groupoid}. 
\end{itemize}
Let us explain what we mean by this. By a \emph{fibered groupoid} we mean a 
groupoid $\G\tto M$ which is fibered over $B$:
\[
\xymatrix{\G\ar@<.5ex>[rr]\ar@<-.5ex>[rr]\ar[dr]& &M\ar[dl]\\ & B}
\]
Here, $\G$ and $B$ are fibered as well as all structure maps.  So a
fibered groupoid maybe thought of as a fiber bundle with fiber type a
groupoid where the structure group acts by groupoid automorphism. Then
by a \emph{fibered symplectic groupoid} we mean a fibered groupoid
$\G$ whose fiber type is a symplectic groupoid. For a fibered
symplectic groupoid, $\G\to B$ is a symplectic fibration and the
symplectic structure is compatible with the groupoid structure. The
base $p:M\to B$ of a fibered symplectic groupoid has a natural
structure of a Poisson fibration. Moreover, the fibers of $\G\to B$
are symplectic groupoids over the fibers of $p:M\to B$, inducing the
Poisson structures on the fibers. In particular, the fibers $p:M\to B$
are integrable Poisson manifolds. Conversely, we will show the
following:

\begin{theorem}
Let $p:M\to B$ be a Poisson fibration with fiber type $(F,\pi)$ an
integrable Poisson manifold. There exists a unique (up to isomorphism)
source 1-connected fibered symplectic groupoid integrating $p:M\to
B$. Moreover, this symplectic fibration always admits a coupling 2-form.
\end{theorem}

It was proved in \cite{BCWZ} that the global objects integrating Dirac 
structures are presymplectic groupoids. Therefore, if $p:M\to B$ is a
Poisson fibration which admits a coupling Dirac structure, there are two
natural groupoids associated with it:
\begin{enumerate}[(i)]
\item The source 1-connected fibered symplectic groupoid integrating 
$p:M\to B$.
\item The source 1-connected presymplectic groupoid integrating the
coupling Dirac structure of $p:M\to B$.
\end{enumerate}
The precise relationship between these two objects is more involved, 
and it would take us too far afield, so it will be discussed elsewhere.

%%%%%%%%%%%%%%%%%%%%%%%%%%%%%%%%%%%%
%%%%%%%%%%%%%%%%%%%%%%%%%%%%%%%%%%%%
%%%%%%%%%%%%%%%%%%%%%%%%%%%%%%%%%%%%
\section{Connections and Dirac structures}%
\label{sec:Dirac:connections}      %
%%%%%%%%%%%%%%%%%%%%%%%%%%%%%%%%%%%%
%%%%%%%%%%%%%%%%%%%%%%%%%%%%%%%%%%%%
%%%%%%%%%%%%%%%%%%%%%%%%%%%%%%%%%%%%
One can define a connection on a fibration by specifying on the total
space an almost Dirac structure of a special type. 
In this section we make a preliminary study of connections 
on fibrations induced by Dirac structures. 

%%%%%%%%%%%%%%%%%%%%%%%%%%%%%%%%%%%%%%%%%%%%%%%%%%%%%%%%%%%%
\subsection{Connections defined by almost Dirac structures}%
\label{sub:sec:connections}                                %
%%%%%%%%%%%%%%%%%%%%%%%%%%%%%%%%%%%%%%%%%%%%%%%%%%%%%%%%%%%%

In what follows, by a \emph{fibration} we always mean a locally trivial 
fiber bundle. Given such a fibration $p:M\to B$, we denote by $F_b$ the fiber
over $b\in B$, and we let $\Vert\subset TM$ be the \emph{vertical sub-bundle} 
whose fibers are $\Vert_x\equiv\Ker\d_x p=T_x F_{p(x)}$. By a 
\emph{connection} $\Gamma$ on $p:M\to B$ we mean an Ehresmann 
connection, i.e., a distribution $\Gamma:x\mapsto \Hor_x$ in $TM$ 
which splits the tangent bundle as
\begin{equation}
\label{eq:connection}
 T_xM=\Hor_x\oplus \Vert_x,
\end{equation}
and satisfies the following lifting property:
\begin{description}
\item[Lifting Property] for each $x_0\in M$ and each curve $\gamma:[0,1]\to B$ starting
at $b_0=p(x_0)$ there exists an integral curve
$\widetilde{\gamma}:[0,1]\to M$ of $\Gamma$ starting at $x_0$ and
covering $\gamma$.
\end{description}

\begin{remark}
Given a $C^1$-path $\gamma:[0,1]\to B$ and $x\in F_{\gamma(0)}$ the
splitting (\ref{eq:connection}) guarantees that there exists a unique
\emph{horizontal lift} $\widetilde{\gamma}_x:[0,\eps]\to M$ starting
at $x$. The Lifting Property says that we can take 
$\eps=1$ (\footnote{Note that our curves are always
parameterized in the interval $[0,1]$.}).
\end{remark}

In the usual way, one obtains the notion of \emph{parallel transport} 
of fibers:  given a piecewise $C^1$-path on the base 
$\gamma:[0,1]\to B$ we have the diffeomorphism:
\[
\phi_{\gamma}(t):F_{\gamma(0)}\to F_{\gamma(t)},\qquad 
x\mapsto \widetilde{\gamma}_x(t). 
\]
When $\gamma$ is a loop, we call $\phi_\gamma(1)$ the \emph{holonomy} 
of $\Gamma$ along $\gamma$. For $b\in B$, the \emph{holonomy group} with base point
$b$ is the subgroup $\Phi(b)\subset\Diff(F_b)$ formed by all
holonomy transformations $\phi_\gamma(1)$, where $\gamma:[0,1]\to B$
is any loop based at $b$. Note that to define the concatenation of loops
we need to reparameterize our curves, but this causes no problem since
``horizontal lift commutes with concatenation'' and, hence, two paths
differing by a reparameterization determine the same holonomy transformation.

Now our basic observation is that connections can be defined by
specifying on the total space of the fibration an almost Dirac 
structure of a special type. Namely:

\begin{definition}
Let $p:M\to B$ be a fibration. An almost Dirac structure $L$ on the
total space of the fibration is called \textbf{fiber non-degenerate} if 
\begin{equation}
\label{eq:fiber:non:degnrt}
(\Vert\oplus\Vert^0)\cap L=\{0\}.
\end{equation}
\end{definition}

In fact, we have:

\begin{prop}
Let $p:M\to B$ be a fibration and $L$ a fiber non-degenerate almost
Dirac structure. Then $L$ defines a connection $\Gamma_L$ with
horizontal space:
\begin{equation}
\label{eq:hor:space}
\Hor=\set{X\in TM:\exists \al\in(\Vert)^0, (X,\al)\in L}.
\end{equation}
Moreover, any connection on $p:M\to B$ can be obtained in this way.
\end{prop}

\begin{proof}
Take a vector space $W$, with a maximal isotropic subspace $L\subset
W\times W^*$, and let $V\subset W$ be a subspace such that:
\[ (V\times V^0)\cap L=\{0\}.\]
We claim that
\[ W=V\oplus H,\]
where $H=\set{w\in W:\exists \xi\in V^0, (w,\xi)\in L}$. If we take
$W=T_x M$, $L=L_x$ and $V=\Vert_x$, this claim yields the first part
of the proposition.

To prove the claim, we start by checking that $V\cap H=\{0\}$. In
fact, if $v\in H$ there exists $\xi\in V^0$ such that $(v,\xi)\in
L$. Hence, if $v\in V\cap H$ we obtain:
\[ (v,\xi)\in (V\times V^0)\cap L=\{0\},\]
so we must have $v=0$. It remains to check that $W=V+H$. First, we
observe that since $L$ and $V\times V^0$ intersect trivially and $\dim
L=\dim (V\times V^0)=\dim W$, we must have:
\[ W\times W^*=(V\times V^0)\oplus L.\]
Therefore, for any $w\in W$ we have a decomposition:
\[ (w,0)=(v,\xi)+(h,\eta),\]
where $v\in V$, $\xi\in V^0$ and $(h,\eta)\in L$. It follows that
$\eta=-\xi\in V^0$, so that $h\in H$. We conclude that $w=v+h$, with
$v\in V$ and $h\in H$, as claimed.

Conversely, if $\Gamma$ is a connection with horizontal distribution 
$\Hor$, then $L=\Hor\oplus\Hor^0$ defines a fiber
non-degenerate almost Dirac structure whose associated connection is
$\Gamma$. Hence, every connection arises in this way.
\end{proof}

%%%%%%%%%%%%%%%%%%%%%%%%%%%%%%%%%%%%%%%%%%%%%%%%%%%%%%%%%%%%%%%%%%%
\subsection{The horizontal 2-form and the vertical bivector field}%
\label{sub:sec:hor:ver:form}                                      %
%%%%%%%%%%%%%%%%%%%%%%%%%%%%%%%%%%%%%%%%%%%%%%%%%%%%%%%%%%%%%%%%%%%

Note that two fiber non-degenerate almost Dirac structures $L_1$ and
$L_2$ may lead to the same connection $\Gamma_{L_1}=\Gamma_{L_2}$. In
fact, as we will see now, there is more structure associated with the
specification of a fiber non-degenerate almost Dirac structure.

\begin{prop}
  Let $p:M\to B$ be a fibration and $L$ a fiber non-degenerate almost
  Dirac structure, with associated connection $\Gamma_L$. Then the horizontal
  distribution $\Hor$ is contained in the characteristic distribution of
  $L$, and the pull-back of the natural 2-form yields a smooth 2-form
  $\omega_L\in\Omega^2(\Hor)$.
\end{prop}

We will refer to $\omega_L$ as the \textbf{horizontal 2-form of $L$}.

\begin{proof}
Fix a fiber non-degenerate almost Dirac structure $L$ on the total
space of a fibration $p:M\to B$. Relations (\ref{eq:fiber:non:degnrt})
and (\ref{eq:hor:space}) together show that, for each $X\in\Hor$, there
exists a unique $\al\in\Vert^0$ such that $(X,\al)\in L$. One can then
define a skew-symmetric bilinear form $\omega:\Hor\times\Hor\to \Rr$ by:
\begin{equation}
\label{eq:2:form}
\omega_L(X_1,X_2):=\frac{1}{2}\left(\al_1(X_2)-\al_2(X_1)\right),
\end{equation}
with $\al_1,\al_2\in\Vert^0$ the unique elements such that
$(X_1,\al_1),(X_2,\al_2)\in L$. Since $L$ is maximal isotropic we
have:
\[ 0=2\langle (X_1,\al_1),(X_2,\al_2)\rangle_+=\al_1(X_2)-\al_2(X_1),\]
so this two form can also be written:
\begin{equation}
\label{eq:2:form:alt}
\omega_L(X_1,X_2)=\al_1(X_2)=-\al_2(X_1).
\end{equation}
In this way we obtain a smooth 2-form $\omega_L\in\Omega^2(\Hor)$. 

From the definition (\ref{eq:hor:space}) of $\Hor$ it is clear that
the horizontal distribution $\Hor$ is contained in the characteristic
distribution of $L$. From (\ref{eq:2:form:alt}), it is clear
that $\omega_L$ is the pull-back of the natural 2-form on the
characteristic distribution of $L$.
\end{proof}

This construction of the horizontal 2-form can be dualized:

\begin{prop}
  Let $p:M\to B$ be a fibration and let $L$ be a fiber non-degenerate almost
  Dirac structure, with associated connection $\Gamma_L$.  For each
  fiber $i:F_x=p^{-1}(x)\hookrightarrow M$ the pull-back almost Dirac
  structure $i^*L$ is well-defined and coincides with the graph of a
  bivector field $\pi_L\in\X^2(\Vert)$. 
\end{prop}

We will refer to $\pi_L$ as the \textbf{vertical bivector field of $L$}.

\begin{proof}
First observe that the annihilator of the horizontal space is:
\begin{equation}
\label{eq:annh:hor:space}
\Hor^0=\set{\al\in T^*M:\exists X\in \Vert, (X,\al)\in L}.
\end{equation}
Relations (\ref{eq:fiber:non:degnrt}) and (\ref{eq:annh:hor:space})
together show that, for each $\al\in\Hor^0$, there 
exists a unique $X\in\Vert$ such that $(X,\al)\in L$. One can then
define a skew-symmetric bilinear form $\pi_L:\Hor^0\times\Hor^0\to \Rr$ by:
\begin{equation}
\label{eq:2:vector}
\pi_L(\al_1,\al_2):=\frac{1}{2}\left(\al_1(X_2)-\al_2(X_1)\right),
\end{equation}
with $X_1,X_2\in\Vert$ the unique elements such that
$(X_1,\al_1),(X_2,\al_2)\in L$. Since $L$ is maximal isotropic we 
have:
\[ 0=2\langle (X_1,\al_1),(X_2,\al_2)\rangle=\al_1(X_2)-\al_2(X_1),\]
the form $\pi_L:\Hor^0\times\Hor^0\to \Rr$ can also be written:
\begin{equation}
\label{eq:2:vector:alt}
\pi_L(\al_1,\al_2)=\al_1(X_2)=-\al_2(X_1).
\end{equation}
Now we remark that the splitting $TM=\Hor\oplus\Vert$ allows us to
identify $\Hor^0=\Vert^*$, so $\pi_L$ becomes a bivector field on the
fibers of $p:M\to B$. 

Let us fix a fiber $i:F_x=p^{-1}(x)\hookrightarrow M$. Then
$TF_x=\Vert$ and we identify $T^*F_x=\Vert^*\simeq\Hor^0$. The pull-back
Dirac structure $i^*L$ is then given by:
\begin{align*}
  i^*L&=\set{(X,\al|_\Vert)\in \Vert\oplus\Vert^*: (X,\al)\in L}\\
      &=\set{(X,\al)\in \Vert\oplus\Hor^0: X=\pi_L(\al,\cdot)}
      =\text{graph}(\pi_L),
\end{align*}
where for the last inequality we have used (\ref{eq:2:vector:alt}).
\end{proof}

Putting together these results, we conclude that:

\begin{cor}
To a fiber non-degenerate almost Dirac structure $L$ on a
fibration $p:M\to B$ there is associated the following data:
\begin{itemize}
\item A connection $\Gamma_L$ on $p:M\to B$.
\item A horizontal 2-form $\omega_L\in\Omega^2(\Hor)$.
\item A vertical bivector field $\pi_L\in\X^2(\Vert)$.
\end{itemize}
Conversely, every such triple $(\Gamma,\omega,\pi)$ on a
fibration $p:M\to B$ determines a unique fiber non-degenerate 
almost Dirac structure $L$, which is given by:
\begin{equation}
  \label{eq:Dirac:components}
  L=\graph(\pi_L)\oplus\graph(\omega_L).
\end{equation}
\end{cor}

%%%%%%%%%%%%%%%%%%%%%%%%%%%%%%%%%%%%%%%%%%%%%%%%%%%%%%%%
\subsection{Fiber non-degenerate Dirac structures}     %
\label{sub:sec:dirac:fibrations}                                       %
%%%%%%%%%%%%%%%%%%%%%%%%%%%%%%%%%%%%%%%%%%%%%%%%%%%%%%%%

The next natural question is: Given a fiber non-degenerate almost 
Dirac structure $L$ on a fibration $p:M\to B$, what are the 
conditions on the associated triple $(\Gamma_L,\omega_L,\pi_L)$ 
that guarantee that $L$ is integrable, i.e., is a Dirac structure? 

Recall (see \cite{Cou}) that the obstruction to integrability for an almost Dirac 
structure $L$ is a 3-form $T_L\in\Omega^3(L)$, which is defined 
on sections $s_1,s_2,s_3\in\Gamma(L)$ by:
\begin{equation}
  \label{eq:T:Dirac}
  T_L(s_1,s_2,s_3):=\langle\llbracket s_1,s_2\rrbracket,s_3\rangle_+
\end{equation}
where:
\begin{itemize}
\item $\llbracket \cdot,\cdot\rrbracket$ denotes the \emph{Courant bracket}, 
on $\X(M)\oplus\Omega^1(M)$, given by:
\begin{equation}
  \label{eq:Courant:bracket}
  \llbracket(X,\al),(Y,\be)\rrbracket:=
  ([X,Y],\Lie_X\be-\Lie_Y\al+\d\langle (X,\al),(Y,\be)\rangle_-)
\end{equation}
\item $\langle\cdot,\cdot\rangle_+$ denotes the \emph{natural pairing} 
on $\X(M)\times\Omega^1(M)$, defined by:
\begin{equation}
 \label{eq:pairings}
\langle (X,\al),(Y,\be)\rangle_+:=\frac{1}{2}\left(i_Y\al+ i_X\be\right).
\end{equation}
\end{itemize}

Our next result gives the 3-form $T_L$ of a fiber non-degenerate almost Dirac structure 
$L$ in terms of the geometric data $(\Gamma_L,\pi_L,\omega_L)$. We need to introduce
some notation. 

For a vector field $v\in\X(B)$ we
denote by $\widetilde{v}\in\X(M)$ its horizontal lift, and we let
$\pi_L^\#:\Vert^*\to\Vert$ denote the bundle map induced by the vertical
bivector field $\pi_L$. Also, we have isomorphisms:
\[ 
\Hor\simeq L\cap(TM\oplus\Vert^0), \quad 
\Hor^0\simeq L\cap(\Vert\oplus T^*M).
\]
These allow us to identify a horizontal vector field $X\in\Gamma(\Hor)$
with a section $s_X=(X,\al)\in \Gamma(L)$, where $\al\in\Gamma(\Vert^0)$, and 
a vertical form $\be\in\Gamma(\Hor^0)$ with a section
$s_\be=(Y,\be)\in\Gamma(L)$, where $Y\in\Gamma(\Vert)$. In this notation we have:

\begin{prop} 
  Let $(\Gamma_L,\pi_L,\omega_L)$ be the geometric data determined by a
  fiber non-degenerate almost Dirac structure $L$ on a fiber bundle
  $p:M\to B$. Then:
  \begin{enumerate}[(i)]
  \item If $\al,\beta,\gamma\in\Gamma(\Hor^0)$ then:
    \[ T_L(s_\al,s_\be,s_\gamma)=\frac{1}{2}[\pi_L,\pi_L](\al,\be,\gamma).\]
   \item If $v\in\X(B)$ and $\beta,\gamma\in\Gamma(\Hor^0)$, then:
    \[ T_L(s_{\widetilde{v}},s_\be,s_\gamma)=
       \frac{1}{2}\Lie_{\widetilde{v}}\pi_L(\be,\gamma).\]
   \item If $v_1,v_2\in\X(B)$ and $\gamma\in\Gamma(\Hor^0)$, then:
    \[ T_L(s_{\widetilde{v}_1},s_{\widetilde{v}_2},s_\gamma)=
       \frac{1}{2}(\gamma(\Omega_{\Gamma_L}(\widetilde{v}_1,\widetilde{v}_2))+
       \pi_L(\d i_{\widetilde{v}_1} i_{\widetilde{v}_2}\omega_L,\gamma))\]
    where $\Omega_{\Gamma_L}$ is the curvature 2-form of $\Gamma_L$.
  \item If $v_1,v_2,v_3\in\X(B)$, then:
    \[ T_L(s_{\widetilde{v}_1},s_{\widetilde{v}_2},s_{\widetilde{v_3}})=
       \frac{1}{2}\d_{\Gamma_L}\omega_L(\widetilde{v}_1,\widetilde{v}_2,\widetilde{v}_3),\]
    where $\d_{\Gamma_L}:\Omega^\bullet(\Hor)\to\Omega^{\bullet+1}(\Hor)$ is
    the differential induced by $\Gamma_L$.
  \end{enumerate}
\end{prop}

The proofs are routine calculations so we omit them. For a different approach to
the geometric data we refer the reader to \cite{DuWa}.

Now observe that sections of the form $s_{\tilde{v}}$ and $s_\al$ with $v\in\X(B)$
and $\al\in\Gamma(\Hor^0)$ generate $\Gamma(L)$, as a
$C^\infty(M)$-module. Therefore, as a corollary, we obtain the
conditions that the triple $(\Gamma_L,\pi_L,\omega_L)$ must satisfy
for the associated $L$ to be a Dirac structure:

\begin{cor} 
\label{cor:obstr:Dirac}
  Let $(\Gamma_L,\pi_L,\omega_L)$ be the geometric data determined by a
  fiber non-degenerate almost Dirac structure $L$ on a fiber bundle
  $p:M\to B$. Then $L$ is Dirac iff the following conditions hold:
  \begin{enumerate}[(i)]
    \item $\pi_L$ is a vertical Poisson structure: $[\pi_L,\pi_L]=0$.
    \item Parallel transport along $\Gamma_L$ preserves the vertical
      Poisson structure $\pi_L$:
      \[ \Lie_{\widetilde{v}}\pi_L=0,\quad \forall v\in\X(B).\]
    \item The horizontal 2-form $\omega_L$ is closed: $\d_\Gamma\omega_L=0$.
    \item The following \textbf{curvature identity} is satisfied:
      \begin{equation}
        \label{eq:curvature:ident}
        \Omega_\Gamma(v_1,v_2)=\pi_L^\#(\d i_{\widetilde{v}_1}
        i_{\widetilde{v}_2}\omega_L),\quad \forall v_1,v_2\in\X(B).
      \end{equation}
    \end{enumerate}
\end{cor}

The curvature identity (\ref{eq:curvature:ident}) expresses the fact
that the curvature 2-form of the connection associated with a fiber
non-degenerate Dirac structure $L$ takes values in the vertical
Hamiltonian vector fields. We will explore this property later in our
study of Poisson fibrations.

%%%%%%%%%%%%%%%%%%%%%%%%%%%%%%%%%%%%%%%%%%%%%%%%%%%%%%%%
\subsection{Presymplectic forms and Poisson structures}%
\label{sub:sec:examples}                               %
%%%%%%%%%%%%%%%%%%%%%%%%%%%%%%%%%%%%%%%%%%%%%%%%%%%%%%%%

Let us illustrate the previous results with the two extreme cases of
Dirac structures determined by Poisson and presymplectic structures. 
\vskip 10 pt

If $L$ is determined by a presymplectic form, one checks easily:

\begin{prop}
Let $\Omega$ be a presymplectic form on the total space of a fibration
$p:M\to B$. Then $L=\graph(\Omega)$ is a fiber non-degenerate Dirac 
structure iff the pull-back of $\Omega$ to each fiber is non-degenerate. 
\end{prop}

In this case, the vertical Poisson structure $\pi_L$ 
is non-degenerate on the fibers and coincides with the inverse of the
restriction of $\Omega$ to the fibers. 

The converse is also true: a fiber non-degenerate Dirac structure $L$
for which the vertical Poisson structure $\pi_L$ is non-degenerate on
the fibers, is determined by a presymplectic form $\Omega$. In fact,
it follows from (\ref{eq:Dirac:components}) that:
\[ \Omega=\omega_L\oplus(\pi_L)^{-1}.\]
Hence, fiber non-degenerate presymplectic forms are presymplectic
forms which restrict to symplectic forms on the fibers.
\vskip 10 pt

Dually, for Poisson structures, it is also immediate to check:

\begin{prop}
Let $\Pi$ denote a Poisson structure on the total space of a fibration
$p:M\to B$. Then $L=\graph(\Pi)$ is a fiber non-degenerate Dirac structure 
iff $\Pi$ is horizontal non-degenerate, i.e., 
$\Pi|_{\Vert^0}:\Vert^0\times\Vert^0\to\Rr$ is a non-degenerate bilinear form.
\end{prop} 

In this case, the horizontal 2-form $\omega_L$ is non-degenerate: in
fact, $\Pi$ gives an isomorphism $\Vert^0\to\Hor$, and under this 
isomorphism $\omega_L$ coincides with the restriction $\Pi|_{\Vert^0}$.

The converse is also true: a fiber non-degenerate Dirac structure $L$
for which the horizontal 2-form $\omega_L$  is non-degenerate, is
a Poisson structure $\Pi$. In fact, it follows from
(\ref{eq:Dirac:components}) that: 
\[ \Pi=(\omega_L)^{-1}\oplus \pi_L.\]
Hence, fiber non-degenerate Poisson structures are the same thing as
the horizontal non-degenerate Poisson structures of Vorobjev (\cite{Vor}).

%%%%%%%%%%%%%%%%%%%%%%%%%%%%%%%%%%%%
%%%%%%%%%%%%%%%%%%%%%%%%%%%%%%%%%%%%
%%%%%%%%%%%%%%%%%%%%%%%%%%%%%%%%%%%%
\section{Poisson fibrations}       %
\label{Poisson:fibrations}         %
%%%%%%%%%%%%%%%%%%%%%%%%%%%%%%%%%%%%
%%%%%%%%%%%%%%%%%%%%%%%%%%%%%%%%%%%%
%%%%%%%%%%%%%%%%%%%%%%%%%%%%%%%%%%%%

In this section we will study Poisson fibrations and their
relationship to fiber non-degenerate Dirac structures.

%%%%%%%%%%%%%%%%%%%%%%%%%%%%%%%%%%%%%%%%%%%%%%%%%%%%%%%%
\subsection{Poisson and symplectic fibrations}         %
\label{sub:sec:fibrations}                             %
%%%%%%%%%%%%%%%%%%%%%%%%%%%%%%%%%%%%%%%%%%%%%%%%%%%%%%%%

Let $(F,\pi)$ be a Poisson manifold. We denote by $\Diff_{\pi}(F)$ the
group of Poisson diffeomorphisms of $F$. This is the subgroup of
$\Diff(F)$ formed by all diffeomorphisms $\phi:F\to F$ such that:
\[ (\phi)_*\pi=\pi.\]
We are interested in the following class of fibrations:

\begin{definition}
A \textbf{Poisson fibration} $p:M\to B$ is a locally trivial fiber
bundle, with fiber type a Poisson manifold $(F,\pi)$ and with
structure group a subgroup $G\subset \Diff_{\pi}(F)$. More precisely, $p$ is a submersion such that there exists a covering $\{U_i\}_i$ of $B$, and trivializations $\phi_i:p^{-1}(U_i)\:\to U_i\times F$, with transition functions ${\phi_j\circ\phi_{i}^{-1}}_{|\{x\}\times F}$, $x\in\mathcal{U}_i\cap\mathcal{U}_j$ belonging to $G$. When $\pi$ is
symplectic the fibration is called a \textbf{symplectic fibration}.
\end{definition}

If $p:M\to B$ is a Poisson fibration modeled on a Poisson manifold
$(F,\pi)$, each fiber $F_b$ carries a natural Poisson structure
$\pi_b$: if $\phi_i:p^{-1}(U_i)\to U_i\times F$ is a local
trivialization, $\pi_b$ is defined by:
\[ \pi_b=(\phi_i(b)^{-1})_*\pi,\]
for $b\in U_i$. It follows from the definition that this 2-vector
field is independent of the choice of trivialization. Note that the Poisson
structures $\pi_b$ on the fibers can be glued to a Poisson structure
$\pi_V$ on the total space of the fibration:
\[ \pi_V(x)=\pi_{p(x)}(x), \quad (x\in M).\]
This 2-vector field is \emph{vertical}: $\pi_V$ takes values in
$\wedge^2\Vert\subset\wedge^2TM$. In this way, the fibers $(F_b,\pi_b)$ 
become Poisson submanifolds of $(M,\pi_V)$.

\begin{example}
\label{ex:neighb:leaf}
An important class of Poisson fibrations is obtained as follows. Take
any Poisson manifold $(P,\Pi)$, fix a \emph{closed} symplectic leaf $B$ of $P$, 
and let $p:M\to B$ be a tubular neighborhood of $B$ in $P$. Each fiber 
carries a natural Poisson structure, namely, the transverse Poisson structure
(\cite{Wein}). These transverse Poisson structures are all Poisson
diffeomorphic and it follows from the Weinstein splitting theorem that
this is a Poisson fibration.

More generally, one can take any Dirac manifold $(P,L)$ and fix a
closed presymplectic leaf $B$ of $P$. If
$p:M\to B$ is a tubular neighborhood of $B$ in $P$, then each fiber
carries a natural Poisson structure, called also the transverse
Poisson structure (\cite{DuWa}). These transverse Poisson structures are
all Poisson diffeomorphic and it follows from the generalization of
the Weinstein splitting theorem in \cite{DuWa} that this is a Poisson fibration.
\end{example}

We saw in the previous section than any fiber non-degenerate Dirac
structure $L$ on the total space of a fibration $p:M\to B$ induces
Poisson structures on the fibers. 
In fact, we have the following:

\begin{prop}
\label{prop:compatible}
Let $p:M\to B$ be a fibration with connected base and compact fibers. 
If $L$ is a fiber non-degenerate Dirac structure, then $p:M\to B$ 
admits the structure of a Poisson fibration such that $\pi_V=\pi_L$.
\end{prop}

\begin{proof}
Corollary \ref{cor:obstr:Dirac} gives the vertical Poisson structure. While trivializations are obtained by parallel transport along paths in $B$, induced by the connection ($F$ being any chosen fiber). Compactness of the fibers here just ensures completeness of horizontal lifts.
\end{proof}

Note that a Dirac structure on the total space of a fibration $p:M\to
B$ for which the fibers are Poisson-Dirac submanifolds may fail to be
a Poisson fibration. In other words, Proposition \ref{prop:compatible}
becomes false if one omits the assumption of fiber
non-degeneracy. This is illustrated by the following simple example.

\begin{example}
Take the fibration $p:\Rr^3\to\Rr$ obtained by projection on the
$x$-axis and the Poisson bracket on $\Rr^3$ defined by:
\[ \{x,y\}=\{x,z\}=0,\quad \{y,z\}=x.\]
Then each fiber is a Poisson-Dirac submanifold: the fiber over $x=0$
has the zero Poisson structure, while the fibers over $x\ne 0$ are
symplectic. Since the fibers are not Poisson diffeomorphic,
$p:\Rr^3\to\Rr$ cannot be a Poisson fibration.
\end{example}

%%%%%%%%%%%%%%%%%%%%%%%%%%%%%%%%%%%%%%%%%%%%%%%%%%%%%%%%
\subsection{Coupling Dirac structures}                 %
\label{sub:sec:coupling}                               %
%%%%%%%%%%%%%%%%%%%%%%%%%%%%%%%%%%%%%%%%%%%%%%%%%%%%%%%%

Motivated by Proposition \ref{prop:compatible} we introduce the
following definition:

\begin{definition}
  If $p:M\to B$ is a Poisson fibration, we will say that a fiber
  non-degenerate Dirac structure $L$ is \textbf{compatible with the
  fibration} if $\pi_L=\pi_V$. In this case, we call $L$ a
  \textbf{coupling Dirac structure}.
\end{definition}

Note that in the special case of a symplectic fibration, by
the results of Section \ref{sub:sec:examples}, a coupling
Dirac structure is necessarily given by a presymplectic form
$\Omega$. In this case, Proposition \ref{prop:compatible} is
well-known (see \cite{MaSa}, Lemma 6.2). Moreover, given a symplectic
fibration $p:M\to B$, there are well-known non-trivial obstructions for
the existence of a coupling 2-form $\omega\in\Omega^2(M)$. Our purpose
now is to determine the corresponding obstructions for a general
Poisson fibration and answer the following question:
\begin{itemize}
\item Given a Poisson fibration $p:M\to B$, is there a
  \emph{coupling} Dirac structure $L$ compatible with the fibration?
\end{itemize}

Let us recall the notion of a Poisson connection:

\begin{definition}
A connection $\Gamma$ on a Poisson fibration $p:M\to B$ is called a
\textbf{Poisson connection} if, for every path $\gamma$, parallel transport
\[ 
\phi_\gamma:(F_{\gamma(0)},\pi_{\gamma(0)})\to(F_{\gamma(1)},\pi_{\gamma(1)})
\]
is a Poisson diffeomorphism. 
\end{definition}

Clearly, a connection $\Gamma$ on a Poisson fibration $p:M\to B$ is
Poisson iff
\[ \Lie_{\widetilde{v}}\pi_V=0,\forall v\in\X(B).\]
Hence, by Corollary \ref{cor:obstr:Dirac} (ii), an obvious necessary
condition for the existence of a coupling Dirac structure is the
existence of a Poisson connection. However, one can show that a
Poisson fibration always admits such a connection. In fact, we have:

\begin{prop}
\label{prop:Poisson:connections}
Let $p:M\to B$ be a Poisson fibration. There exists a fiber
non-degenerate almost Dirac structure $L$ on $p:M\to B$ such that:
\begin{enumerate}[(a)]
\item The vertical bivector field $\pi_L$ coincides with $\pi_V$.
\item The connection $\Gamma_L$ is a Poisson connection.
\end{enumerate}
\end{prop}

\begin{proof}
Let $p:M\to B$ be a Poisson fibration with fiber $(F,\pi)$ and choose
choose local trivializations $\phi_i:p^{-1}(U_i)\to U_i\times F$. 
Let $L_i$ be the Dirac structure on $U_i\times F$ obtained by
pull-back of $L_\pi=\graph(\pi)$ under the projection $U_i\times F\to F$:
\[ L_i:=\{((v,w),(0,\eta))\in T(U_i\times F)\oplus T^*(U_i\times F):
w=\pi(\eta,\cdot)\}.\]
Observe that to $L_i$ it is associated the geometric data
$(\Gamma_{L_i},\pi_{L_i},\omega_{L_i})$ where $\Gamma_{L_i}$ is the
canonical flat connection on $U_i\times F\to U_i$, $\pi_{L_i}=\pi$ and
$\omega_{L_i}=0$. Hence $L_i$ is fiber non-degenerate, induces $\pi$
on the fibers, and $\Gamma_{L_i}$ is a Poisson connection.

Next, we choose a partition of unity $\rho_i:B\to\Rr$ subordinated to
the cover $\{U_i\}$, and we define
\[ L:=\sum_i (\rho_i\circ p)\phi_i^*L_i.\]
By this we mean that the associated geometric data
$(\Gamma_L,\pi_L,\omega_L)$ has connection $\Gamma_L=\sum_i
(\rho_i\circ p)\phi_i^*\Gamma_{L_i}$, vertical bivector field 
$\pi_L=\sum_i(\rho_i\circ p)\phi_i^*\pi_{L_i}$ and horizontal 2-form
$\omega_L=\sum_i(\rho_i\circ p)\phi_i^*\omega_{L_i}=0$. 

It is clear that $L$ is a fiber non-degenerate almost Dirac
structure. Moreover, $\pi_L$ coincides with $\pi_V$. Finally, given
$v\in\X(B)$, if $\widetilde{v}_i$ denotes the horizontal lift relative to
$\phi_i^*\Gamma_{L_i}$, we have:
\[ \Lie_{\widetilde{v}_i}\pi_V=0.\]
The horizontal lift of $v$ relative to $\Gamma_L$ is 
$\widetilde{v}=\sum_i(\rho_i\circ p)\widetilde{v}_i$. It follows that:
\begin{align*}
\Lie_{\widetilde{v}}\pi_V&=[\widetilde{v},\pi_V]\\
&=\sum_i[(\rho_i\circ p)\widetilde{v}_i,\pi_V]\\
&=\sum_i\left((\rho_i\circ p)[\widetilde{v}_i,\pi_V]+
\pi_V^\#(\d (\rho_i\circ p))\right)\\
&=\sum_i(\rho_i\circ p)\Lie_{\widetilde{v}_i}\pi_V=0,
\end{align*}
where we have used the fact that $\d (\rho_i\circ p)\in\Vert^0$.
Hence $L$ is the desired almost Dirac structure.
\end{proof}

There is also a procedure to construct coupling Dirac structures due to
A.~Wade (\cite{Wa}), which is entirely analogous to a construction in 
symplectic geometry due to A.~Weinstein (see \cite[Theorem 6.17]{MaSa}):

\begin{theorem}
\label{thm:Yang:Mills}
Let $G\times F\to F$ be a Hamiltonian action of a compact Lie group $G$ on
the Poisson manifold $(F,\pi)$. Every connection on a principal
$G$-bundle $P\to B$ determines a coupling Dirac structure $L$ on the
associated Poisson fibration $P\times_G F\to B$.
\end{theorem}

\begin{proof}
The connection $\Gamma$ on $P$ determines a projection $TP\to\Vert$, along the
horizontal distribution $\Hor$. Hence, there is an injection:
\[ i_\Gamma:\Vert^*\hookrightarrow T^*P,\]
where $\Vert^*$ is the vertical cotangent bundle with typical fiber
$T^*P_b$. Since $\Gamma$ is $G$-invariant, the inclusion is $G$-equivariant.
Therefore, the canonical symplectic form $\omega_{\text{can}}$ in
$T^*P$ induces a closed 2-form on $\Vert^*$ 
\[ \omega_\Gamma=i_\Gamma^*\omega_{\text{can}},\] 
which is $G$-invariant and restricts to the canonical symplectic form
on the fibers $T^*P_b$. Also, the action of $G$ on $\Vert^*$ has a
moment map $\mu_P\circ i_\Gamma:\Vert^*\to\gg^*$, where $\mu_P$ is the
moment map of the lifted cotangent action $G\times T^*P\to T^*P$.

Now let $G\times F\to F$ be a Hamiltonian action with moment map
$\mu_F:F\to\gg^*$. This determines a Hamiltonian $G$-action on the
Dirac manifold $M:=\Vert^*\times F$, with Dirac structure 
\[ L=\graph(\omega_\Gamma)\oplus\graph(\pi).\]
and with moment map 
\[ \mu_M=(\mu_P\circ i_\Gamma)\oplus\mu_F.\]
This action is free and $0$ is a regular value of $\mu_M$. Therefore,
the reduced space 
\[ M_\text{red}:=\mu_M^{-1}(0)/G\simeq P\times_G F,\]
carries a Dirac structure. It is easy to check that this Dirac
structure is fiber non-degenerate and restricts to the canonical
Poisson structures on the fibers. Hence, $L$ is the desired coupling
Dirac structure.
\end{proof}

\begin{remark}
\label{rem:Yang:Mills}
Note that the resulting coupling Dirac structure is presymplectic iff
the fiber $(F,\pi)$ is symplectic. Also, it is easy to check that the
remaining geometric data associated with the coupling Dirac structure
$L$ in the theorem above is the following:
\begin{itemize}
\item The connection $\Gamma_L$ is just the connection on the fiber
  bundle $P\times_G F$ induced from the connection $\Gamma$ on $P$.
\item The horizontal 2-form $\omega_L$ is given by the curvature of
  the connection composed with the moment map
  \[ \omega_L(\tilde{v}_1,\tilde{v}_2)([u,x])=\langle
  \mu_F(x),F_\Gamma(\tilde{v}_1,\tilde{v}_2)_u\rangle,\quad
  ([u,x]\in P\times_G F).\]
\end{itemize}
Hence, the resulting coupling Dirac structure is Poisson iff the
curvature 2-form of the connection is non-degenerate on the image of
$\mu_F$. Such a connection is sometimes called a \emph{fat} connection.
The proof of A.~Wade in \cite{Wa} consists in proving that this data 
satifies the conditions of Corollary \ref{cor:obstr:Dirac} and so 
defines a coupling. 
\end{remark}

%%%%%%%%%%%%%%%%%%%%%%%%%%%%%%%%%%%%%%%%%%%%%%%%%%%%%%
\subsection{Obstruction to the existence of coupling}%
\label{sub:sec:obstr:coupling}                       %
%%%%%%%%%%%%%%%%%%%%%%%%%%%%%%%%%%%%%%%%%%%%%%%%%%%%%%

The construction of Theorem \ref{thm:Yang:Mills} can be extended for
general Poisson fibrations, leading to a characterization of those
fibrations which admit a coupling Dirac structure.

Let $L$ be a fiber non-degenerate almost Dirac structure on a Poisson
fibration $p:M\to B$. We will say that $L$ is compatible with the
fibration if $\pi_L=\pi_V$. By Proposition \ref{prop:Poisson:connections}, 
any Poisson fibration admits a compatible $L$ such that $\Gamma_L$ is
a Poisson connection. For $L$ to be Dirac this connection must have a much
more constrained holonomy:

\begin{theorem}
\label{thm:coupling}
Let $p:M\to B$ be a Poisson fibration and let $L$ be a compatible
almost Dirac structure such that $\Gamma_L$ is a Poisson
connection. Then the following statements are equivalent:
\begin{enumerate}[(i)]
\item $L$ is a Dirac structure: $T_L=0$.
\item For every base point $b\in B$, the action of the holonomy group
  $\Phi(b)$ of $\Gamma_L$ on the fiber $F_b$ is Hamiltonian.
\end{enumerate}
\end{theorem}

A proof of this result will be given in the next paragraph, using a
Poisson gauge theory. An immediate corollary is the following result
(see, also, \cite[Theorem 3.4]{Wa}):

\begin{cor}\label{cor:coupling}
Let $(F,\pi)$ be a compact Poisson manifold whose first Poisson cohomology
group vanishes: $H^1_\pi(F)=0$. Then any Poisson fibration $p:M\to B$
with fiber $(F,\pi)$, and finite dimensional structure group admits a 
coupling Dirac structure.
\end{cor}

\begin{proof}
Since $H^1_{\pi}(F)=0$, the same holds for the fibers $F_b$, and it follows
that any Poisson action on the fibers is Hamiltonian.  Now apply
Proposition \ref{prop:Poisson:connections} to construct a fiber
non-degenerate almost Dirac structure $L$ compatible with the
fibration and such that $\Gamma_L$ is a Poisson connection. Using the
implication (ii) $\Rightarrow$ (i) in Theorem \ref{thm:coupling}, we
conclude that $L$ is a coupling Dirac structure. Compactness of $F$ ensures 
completeness of the connection.
\end{proof}

Notice that the condition that the holonomy group $\Phi(b)$ of
$\Gamma_L$ acts in a Hamiltonian fashion on the fiber $F_b$ is a
property of its connected component of the identity $\Phi(b)^0$. This
connected component is known as the \emph{restricted holonomy group}
and is formed by the holonomy homomorphisms $\phi_\gamma$, where
$\gamma$ is a contractible loop based at $b$. 

In particular, each $\phi_\gamma$ with $\gamma$ a contractible loop
based at $b$, lies in the group of Hamiltonian diffeomorphisms
$\Ham(F_b,\pi_b)$, which is known to be a normal subgroup of the group
of Poisson diffeomorphisms $\Diff(F_b,\pi_b)$. The quotient group 
$\Diff(F_b,\pi_b)/\Ham(F_b,\pi_b)$ is known as the group of
\emph{outer Poisson diffeomorphisms}. We conclude that a coupling
Dirac structure $L$ for a Poisson fibration $p:M\to B$ has an
associated \emph{coupling holonomy homomorphism}:
\[ \phi:\pi_1(B,b)\to \Diff(F_b,\pi_b)/\Ham(F_b,\pi_b),\quad
   [\gamma]\mapsto [\phi_\gamma].\] 

\begin{example}
A tubular neighborhood $p:M\to B$ of a symplectic leaf $B$ of a
Poisson manifold $(P,\Pi)$ (see Example \ref{ex:neighb:leaf}) admits
$L_\Pi$ as a coupling Dirac structure. It follows that the connection
$\Gamma_\Pi$ is Poisson and has Hamiltonian holonomy around any
contractible loop in $B$. This can also be proved directly using the
Weinstein splitting theorem.

In general, the holonomy around a non-contractible loop will not be
Hamiltonian and we will have a nontrivial homomorphism 
\[ \phi:\pi_1(B)\to \Diff(F,\pi)/\Ham(F,\pi).\]
This is precisely the (reduced) Poisson holonomy of the leaf $B$ 
introduced in \cite{Fe}.
\end{example}

%%%%%%%%%%%%%%%%%%%%%%%%%%%%%%%%%%%%
\subsection{Poisson Gauge Theory}  %
\label{sub:sec:gauge:theory}       %
%%%%%%%%%%%%%%%%%%%%%%%%%%%%%%%%%%%%

We now turn to the proof of Theorem \ref{thm:coupling}. The idea will
be to give an analogue of Theorem \ref{thm:Yang:Mills}, but where the
structure group is allowed to be infinite dimensional. 

We consider a Poisson fibration $p:M\to B$ with fiber type a Poisson
manifold $(F,\pi)$. The structure group of this fibration is the group
$G=\Diff(F,\pi)$ of Poisson diffeomorphisms. The corresponding
principal $G$-bundle is the \textbf{Poisson frame bundle}:
\[ P\to B\]
whose fiber over a point $b\in B$ is formed by all Poisson
diffeomorphisms $u:F\to F_b$. The group $G$ acts on (the right of) $P$ by
pre-composition: 
\[ P\times G\to P:~(u,g)\mapsto u \circ g.\]
Then our original Poisson fiber bundle is canonically isomorphic to
the the associated fiber bundle: $M=P\times_G F$.

Every Poisson connection $\Gamma$ on the Poisson fiber bundle $p:M\to
B$ is induced by a principal bundle connection on $P\to B$. To see
this, observe that the tangent space $T_u P\subset C^\infty(u^*TM)$ at
a point $u\in P$ is formed by the vector fields along $u$, $X(x)\in
T_{u(x)}M$ such that:
\[ d_{u(x)}p\cdot X(x)=\text{constant},\quad \Lie_X\pi_V=0.\]
The Lie algebra $\gg$ of $G$ is the space of Poisson vector fields:
$\gg=\X(F,\pi)$. The infinitesimal action on $P$ is given by:
\[ \rho:\gg\to\X(P),\quad \rho(X)_u= \d u\cdot X,\]
so the vertical space of $P$ is:
\[ \Vert_u=\set{\d u\cdot X: X\in \X(F,\pi)}.\]
Now a Poisson connection $\Gamma$ on $p:M\to B$ determines a
connection in $P\to M$ whose horizontal space is:
\[ \Hor_u=\set{\widetilde{v}\circ u: v\in T_bB},\]
where $u:F\to F_b$ and $\widetilde{v}:F_b\to T_{F_b }M$ denotes the
horizontal lift of $v$.  Clearly, this defines a principal bundle
connection on $P\to B$, whose induced connection on the associated
bundle $M=P\times_G F$ is the original Poisson connection $\Gamma$.

Fix a Poisson connection $\Gamma$ on the Poisson fiber bundle $p:M\to
B$. Recall that the holonomy group $\Phi(b)$ with base point $b\in B$
is the group of holonomy transformations $\phi_\gamma:F_b\to F_b$,
where $\gamma$ is a loop based at $b$. Clearly, we have
$\Phi(b)\subset\Diff(F_b,\pi_b)$. On the other hand, for $u\in P$ we
have the holonomy group $\Phi(u)\subset G=\Diff(F,\pi)$ of the
corresponding connection in $P$ which induces $\Gamma$: it consist of
all elements $g\in G$ such that $u$ and $ug$ can be joined by a
horizontal curve in $P$. Obviously, these two groups are isomorphic, for if
$u:F\to F_b$ then:
\[ \Phi(u)\to \Phi(b),\ g\mapsto u\circ g\circ u^{-1},\]
is an isomorphism.

The curvature of a principal bundle connection is a $\gg$-valued
2-form $F_\Gamma$ on $P$ which transforms as:
\[ R_g^*F_\Gamma=\Ad(g^{-1})\cdot F_\Gamma,\quad (g\in G).\]
Therefore, we can also think of the curvature as a 2-form $\Omega_L$
with values in the adjoint bundle $\gg_P:=P\times_G\gg$. In the case of the
Poisson frame bundle, the adjoint bundle has fiber over $b$ the space
$\X(F_b,\pi_b)$ of Poisson vector fields on the fiber. Hence the
curvature of our Poisson connection can be seen as a 2-form
$\Omega_\Gamma:T_bB\times T_bB\to\X(F_b,\pi_b)$. The two curvature
connections are related by:
\begin{equation}
\label{eq:curv:relations}
\Omega_\Gamma=\d u\circ F_\Gamma\circ u^{-1}.
\end{equation}
Finally, it is easy to check that, in fact, we have:
\[
\Omega_\Gamma(v_1,v_2)=[\widetilde{v}_1,\widetilde{v}_2]-\widetilde{[v_1,v_2]},
\]
which is the expression we have used before for the curvature. 

After these preliminarities, we can now proceed to the proof.

\begin{proof}[Proof of Theorem \ref{thm:coupling}] We will prove the
  two implications separately.

(i) $\Rightarrow$ (ii). Let us start by observing that given any $u\in P$,
a Poisson diffeomorphism $u:F\to F_b$, the curvature identity
(\ref{eq:curvature:ident}) together with (\ref{eq:curv:relations})
shows that, for any $v_1,v_2\in T_bB$, the vector field
$F_\Gamma(v_1,v_2)_u\in\gg=\X(F,\pi)$ is Hamiltonian:
\[ 
F_\Gamma(v_1,v_2)_u=
\pi^{\#}\d (\omega_L(\widetilde{v}_1,\widetilde{v}_2)\circ u).
\]

Now fix $u_0\in P$. The Holonomy Theorem states that the Lie algebra
of the holonomy group $\Phi(u_0)$ is generated by all values
$F_\Gamma(v_1,v_2)_u$, with $u\in P$ any point that can be connected
to $u_0$ by a horizontal curve. Hence, we can define a moment map
$\mu_F:F\to (\text{Lie}(\Phi(u_0))^*)$ for the action of $\Phi(u_0)$
on $F$ by:
\begin{equation}
\label{eq:moment:map}
\langle \mu_F(x),F_\Gamma(v_1,v_2)_u\rangle
=\omega_L(\widetilde{v}_1,\widetilde{v}_2)_{u(x)}.
\end{equation}
This shows that the action of $\Phi(u_0)$ on $(F,\pi)$ is
Hamiltonian, and so (ii) holds (recall the comments above about the
relationship between the holonomy groups $\Phi(b)$ and $\Phi(u)$).
\vskip 10 pt

(ii) $\Rightarrow$ (i). Again we fix $u_0\in P$, and we assume now
that the action of $\Phi(u_0)$ on $(F,\pi)$ is Hamiltonian with moment
map $\mu_F:F\to (\text{Lie}(\Phi(u_0))^*)$. By the Reduction Theorem
we can reduce the principal Poisson frame bundle to a principal
$\Phi(u_0)$-bundle $P'\to B$. Now we can apply (the infinite
dimensional version) of Theorem \ref{thm:Yang:Mills} to produce a
coupling Dirac structure on the associated Poisson fiber bundle
$p:M\to B$. Instead, if the reader does not like an infinite
dimensional argument, he can check by himself that the geometric data
formed by the connection $\Gamma_L$, the vertical Poisson vector field
$\pi_V$ and the 2-form $\omega_L$ defined from (\ref{eq:moment:map})
(we are now given $\mu_F$ and we define $\omega_L$) satisfy the
conditions of Corollary \ref{cor:obstr:Dirac}.
\end{proof}

\begin{remark}
Note that our proof really shows that for \emph{any} Poisson fibration 
the coupling Dirac structure arises as in the construction of Theorem
\ref{thm:Yang:Mills}. Relation (\ref{eq:moment:map}) between the
moment map $\mu_F$, the curvature of the connection, and the
horizontal 2-form $\omega_L$, was already present there (see Remark
\ref{rem:Yang:Mills}).
\end{remark}
%%%%%%%%%%%%%%%%%%%%%%%%%%%%%%%%%%%%%%%%%%%%
%%%%%%%%%%%%%%%%%%%%%%%%%%%%%%%%%%%%%%%%%%%%
%%%%%%%%%%%%%%%%%%%%%%%%%%%%%%%%%%%%%%%%%%%%
\section{Integration of Poisson fibrations}%
\label{sec:integ:fibrations}               %
%%%%%%%%%%%%%%%%%%%%%%%%%%%%%%%%%%%%%%%%%%%%
%%%%%%%%%%%%%%%%%%%%%%%%%%%%%%%%%%%%%%%%%%%%
%%%%%%%%%%%%%%%%%%%%%%%%%%%%%%%%%%%%%%%%%%%%

In this section, we study the integration of Poisson
fibrations. Just as Poisson manifolds integrate to symplectic
groupoids, we will see that Poisson fibrations integrate to fibered
symplectic groupoids.

%%%%%%%%%%%%%%%%%%%%%%%%%%%%%%%%%%%%%%%%%%%
\subsection{Fibered symplectic groupoids} %
\label{sub:sec:fibered:groupoids}         %
%%%%%%%%%%%%%%%%%%%%%%%%%%%%%%%%%%%%%%%%%%%

Recall that for us a fibration always means a locally trivial fiber
bundle. If we fix a base $B$, we have a category $\mathbf{Fib}$ of fibrations over
$B$, where the objects are the fibrations $p:M\to B$ and the morphisms
are the fiber preserving maps over the identity:
\[
\xymatrix{M_1\ar[rr]^{\phi}\ar[rd]_{p_1}&&M_2\ar[ld]^{p_2}\\ &B&}
\]

A \emph{fibered groupoid} is an internal groupoid in $\mathbf{Fib}$, i.e.,
an internal category where every morphism is an isomorphism. This means
that both the total space $\G$ and the base $M$ of a fibered groupoid are
fibrations over $B$ and all structure maps are fibered maps. For example, the source
and target maps are fiber preserving maps over the identity:
\[
\xymatrix{\G\ar@<.5ex>[rr]\ar@<-.5ex>[rr]\ar[dr]& &M\ar[dl]\\ & B}
\]
In particular, each fiber of $\G\to M$ is a groupoid over a fiber of
$M\to B$. Moreover, the orbits of $\G$ lie inside the fibers of the
base $M$.

A general procedure to construct fibered Lie groupoids is as
follows. Let $P\to B$ be a principal $G$-bundle and assume that $G$
acts on a groupoid $\F\tto F$ by groupoid automorphisms. Then the
associated fiber bundles $\G=P\times_G\F$ and $M=P\times_G F$ are the
spaces of arrows and objects of a fibered Lie groupoid. Clearly, every
fibered Lie groupoid is of this form provided we allow infinite
dimensional structure groups. We will say that the fibered Lie groupoid
$\G$ has \emph{fiber type} the Lie groupoid $\F$.

\begin{definition}
A \textbf{fibered symplectic groupoid} is a fibered Lie groupoid $\G$
whose fiber type is a symplectic groupoid $(\F,\omega)$.
\end{definition}

Therefore, if $\G$ is a fibered symplectic groupoid over $B$, then
$\G\to B$ is a symplectic fibration, and each symplectic fiber $\F_b$
is in fact a symplectic groupoid over the corresponding fiber $F_b$ of
$M\to B$.

\begin{prop}
\label{prop:Lie:I}
The base $M\to B$ of a fibered symplectic groupoid $\G\tto M$ has a natural
structure of a Poisson fibration.
\end{prop}

\begin{proof}
Note that (i) the base of any symplectic groupoid has a natural
Poisson structure for which the source (respectively, the target) is a
Poisson (respectively, anti-Poisson) map, and (ii) any symplectic
groupoid isomorphism between two symplectic groupoids covers a Poisson
diffeomorphism of the base Poisson manifolds. Hence, each fiber of
the base $M\to B$ of a fibered symplectic groupoid carries a natural
Poisson structure, and a trivialization of the fibered symplectic
groupoid covers a trivialization of $M\to B$ whose transition functions are
Poisson diffeomorphisms of the fibers. Therefore the result follows.
\end{proof}

Note that the fibers $p:M\to B$ are integrable Poisson
manifolds.

%%%%%%%%%%%%%%%%%%%%%%%%%%%%%%%%%%%%%%%%%%%%%%%
\subsection{Integration of Poisson fibrations}% %
\label{sub:sec:integr:fibrations}             %
%%%%%%%%%%%%%%%%%%%%%%%%%%%%%%%%%%%%%%%%%%%%%%%

We just saw that the base of a fibered symplectic groupoid is a Poisson fibration.
We will say that a fibered symplectic groupoid $\G\to B$ \emph{integrates} a Poisson
fibration $p:M\to B$ whenever this fibration is (Poisson) isomorphic
to the Poisson fibration determined by $\G\to B$. If such a fibered symplectic
groupoid exists we say that the Poisson fibration is \emph{integrable}. Note that the
fiber type $\F$ of $\G$ is a symplectic groupoid integrating the fiber type
$(F,\pi)$ of $p:M\to B$.

\begin{theorem}
\label{thm:Lie:III}
A Poisson fibration is integrable iff its fiber type is an integrable Poisson manifold.
There exists a 1:1 correspondence between, source 1-connected, fibered symplectic groupoids
and integrable Poisson fibrations.
\end{theorem}

In one direction, the proof follows from Proposition \ref{prop:Lie:I}.
In the other direction, we will offer two proofs. The first proof is an heuristic
proof that uses Poisson gauge theory. The second proof uses the approach
to integrability through cotangent paths developed in \cite{CrFe2,CrFe1}.

\begin{proof}[Heuristic proof via gauge theory]
Given a Poisson fibration $p:M\to B$ with fiber type $(F,\pi)$, we start
by writing it as an associated fiber bundle:
\[
M=P\times_G F
\]
where $P$ is the Poisson frame bundle and $G\subset\Diff(F,\pi)$ is the
structure group of the fibration.

Since $(F,\pi)$ is integrable, there exists a unique source
1-connected symplectic groupoid $\F\tto F$ which integrates $(F,\pi)$.
The action of $G$ on $F$ lifts to an action of $G$ on
$\F$ by symplectic groupoid automorphisms (see \cite{FeOrRa}). Hence, we
can form the associated bundle:
\[ \G=P\times_G\F.\]
Since the action of $G$ on $\F$ is by groupoid automorphisms, $\G\to B$ becomes
a fibered groupoid over $M\to B$. Since this action is by symplectomorphisms,
$\G\to B$ becomes a symplectic fibration. Since the groupoid structures and the
symplectic structure on the fibers are compatible, $\G$ is a source 1-connected
fibered symplectic groupoid with fiber type $\F$. It should be clear that the Poisson
fibration determined by $\G$ is isomorphic to the original fibration.
\end{proof}

\begin{remark}
Note that this heuristic proof becomes a real proof if the structure group $G$
of the fibration is a finite dimensional Lie group. We will illustrate this
below in Example \ref{sub:sec:integr:exampl}.
\end{remark}

\begin{proof}[Proof of Theorem \ref{thm:Lie:III}]
Given a Poisson fibration $p:M\to B$ with fiber type $(F,\pi)$, we denote
by $\Sigma(M)\tto M$ the symplectic groupoid that integrates
the vertical Poisson structure $\pi_V$. Note that since we assume that
$(F,\pi)$ is integrable, we have that $(M,\pi_V)$ is integrable, so
that $\Sigma(M)$ is a Lie groupoid.

Let us recall (see \cite{CrFe2,CrFe1} for details and notations) that
$\Sigma(M)$ is the space of equivalence classes of cotangent paths:
\[ \Sigma(M)=\frac{\{a:[0,1]\to T^*M:\pi_V^\sharp(a(t))=\frac{\d}{\d t}p(a(t))\}}
{\{\text{cotangent homotopies}\}} ;\]
where $p:T^*M\to M$ is the cotangent bundle projection. Now we observe that
$\Vert^0\subset T^*M$ is a Lie subalgebroid, which is
in fact, a bundle of Abelian Lie algebras. This is a direct consequence of $\pi_V$
being vertical. Hence, the equivalence classes of cotangent paths with
image in $\Vert^0$ form a closed Lie subgroupoid $\K\subset\Sigma(M)$, which is in fact
a bundle of Abelian Lie groups.

Let us consider the quotient Lie groupoid:
\[ \G:=\Sigma(M)/\K. \]
Notice that $\G$ is fibered over $B$, where the fiber over $b$ is the
symplectic groupoid $\Sigma(F_b)$ integrating the fiber $(F_b,\pi_b)$.
It is easy to check that $\G$ is, in fact, the desired source 1-connected,
fibered symplectic groupoid integrating $p:M\to B$. We leave the details to the reader.
\end{proof}

As we have seen in Proposition \ref{prop:Poisson:connections}, a Poisson fibration
$p:M\to B$ always admits Poisson connections. What does the specification of a
Poisson connection on $p:M\to B$ amounts to in the corresponding fibered symplectic
groupoid $\G$?

\begin{prop}
Let $M\to B$ be a Poisson fibration which integrates to a
source 1-connected fibered symplectic groupoid $\G\to B$. The choice
of a Poisson connection on the Poisson fibration $M\to B$ determines
a coupling form $\Omega$ on the fibered symplectic groupoid $\G\to B$,
and conversely.
\end{prop}

\begin{proof}
We will give a ``gauge theoretical'' proof, which is valid at least in the
case where the structure group is a finite dimensional Lie group. One can also 
give a longer proof using paths, which avoids this assumption.

Hence, assume that:
\[ M=P\times_G F,\]
where $P$ is a principal $G$-bundle and $F$ is a Poisson $G$-space.
As we have mentioned above, the Poisson action $G\times F\to F$ lifts
to an action $G\times \F\to \F$ by automorphisms of the symplectic
groupoid $\F=\Sigma(F)$, which is Hamiltonian with equivariant moment
map $J:\F\to\gg^*$, which is a groupoid cocycle. We have that:
\[ \G=P\times_G \F. \]
Now, Poisson connections $\Gamma$ on $M\to B$ are in 1:1 correspondence
with principal bundle connections on $P$.

To complete the proof we observe that, since the action $G\times \F\to \F$
is Hamiltonian, a choice of a principal bundle connection on $P$ determines
a coupling form on $\G$ and conversely (see Theorem \ref{thm:Yang:Mills}
or \cite[Chapter 6]{MaSa} for more details).
\end{proof}

The next natural question is: what does a coupling Dirac structure on
the Poisson fibration $p:M\to B$ amounts to in the corresponding fibered
symplectic groupoid $\G$? This question is more delicate, and it is intimately
related with the pre-symplectic groupoids integrating Dirac structures described
in \cite{BCWZ}. This will be discussed elsewhere.

%%%%%%%%%%%%%%%%%%%%%%%%%%%%%%%%%%%%%%%%%%%%%%%
\subsection{An Example} 											%
\label{sub:sec:integr:exampl}                 %
%%%%%%%%%%%%%%%%%%%%%%%%%%%%%%%%%%%%%%%%%%%%%%%

Let us denote by $\Ss^3\to\Ss^2$ the Hopf fibration which we view as
a principal $\Ss^1$-bundle $P\to\Ss^2$. We will consider as fiber types
$(F,\pi)$ the following two Poisson $\Ss^1$-manifolds:
\begin{enumerate}[1)]
\item The manifold $F=\Ss^2$, with the standard area form and the
$\Ss^1$-action by rotations around the north-south poles axis;
\item The manifold $F=\su(2)^*\simeq\Rr^3$, with its canonical linear
Poisson structure and the $\Ss^1$-action by rotations around the
$z$-axis;
\end{enumerate}
The corresponding Poisson fibrations $M=P\times_{\Ss^1}F$ are:
\begin{enumerate}[1)]
\item the non-trivial $\Ss^2$-bundle $p:M\to\Ss^2$ (a symplectic fibration), and
\item the non-trivial rank 3 vector bundle $p:E\to\Ss^2$ (a Poisson
fibration which is not symplectic).
\end{enumerate}
The symplectic leaves of $\su(2)^*$ are the concentric spheres around the
origin and the origin itself. Hence the Poisson fibration $p:E\to\Ss^2$ is
foliated by symplectic fibrations isomorphic to $p:M\to\Ss^2$ and the
zero section. Since $\Ss^2$ is symplectic, we have:
\[ H^1_\pi(\Ss^2)\simeq H^1(\Ss^2)=\{0\}.\]
Since $\su(2)$ is semisimple of compact type, we also have:
\[ H^1_\pi(\su(2)^*)=\{0\}.\]
It follows from Corollary \ref{cor:coupling} that both $p:E\to\Ss^2$ and
$p:M\to\Ss^2$ admit coupling Dirac structures. Of course, since $p:M\to\Ss^2$
is a symplectic fibration, its Dirac coupling is actually associated with
a closed 2-form. We let the reader check that the presymplectic leaves of
the Dirac coupling for $p:E\to\Ss^2$ are the symplectic fibrations isomorphic
to $p:M\to\Ss^2$ (with their coupling forms) and the zero section (with the zero
2-form).

Let us now turn to the fibered symplectic groupoids integrating these fibrations.
For that, we use the method in the heuristic proof of Theorem \ref{thm:Lie:III}.
Since the structure group is $\Ss^1$, a finite dimensional Lie group, this is
allowed. We need the source 1-connected symplectic groupoid $\F=\Sigma(F)$
integrating the fiber type, and this is well-known in both examples:
\begin{enumerate}[1)]
\item Since $\Ss^2$ is symplectic and 1-connected, the associated
source 1-connected symplectic groupoid is the pair groupoid
$\Sigma(\Ss^2)=\Ss^2\times\overline{\Ss^2}$, where the bar over the
second factor means that we change the sign of symplectic form.
\item From general facts about linear Poisson structures, the
symplectic groupoid of $\su(2)^*$ is $\Sigma(\su(2)^*)=T^*\SU(2)$,
furnished with the canonical cotangent bundle symplectic structure.
This groupoid is isomorphic to the action groupoid $\SU(2)\ltimes\su(2)^*$
for the coadjoint action of $\SU(2)$ on $\su(2)^*$.
\end{enumerate}

Now we can describe, in both cases, the fibered symplectic groupoid $\G\to\Ss^2$
given by Theorem \ref{thm:Lie:III}.

For the non-trivial $\Ss^2$-fibration $M\to\Ss^2$, the action of $\Ss^1$ on $\Ss^2$
lifts to the diagonal $\Ss^1$-action on $\Ss^2\times\overline{\Ss^2}$, and we
have:
\[ \G(M)=P\times_{\Ss^1} (\Ss^2\times\overline{\Ss^2}),\]
which is a non-trivial symplectic $(\Ss^2\times\overline{\Ss^2})$-fibration over $\Ss^2$
and a groupoid over the non-trivial $\Ss^2$-fibration.

For the rank 3 vector bundle $E\to\Ss^2$, the action of $\Ss^1$ on $\su(2)^*$
lifts to an action on $\SU(2)\times\su(2)^*$, which is trivial on the first factor,
and we have:
\[ \G(E)=P\times_{\Ss^1} (\SU(2)\times\su(2)^*)\simeq E\times \SU(2).\]
Note that, contrary to the case of the Poisson fibrations, the symplectic
fibered groupoid $\G(M)$ does not sit naturally in $\G(E)$. This is because
a Poisson submanifold does not always integrate to a symplectic subgroupoid,
and this is exactly the case with the spheres in $\su(2)^*$.

\bibliographystyle{amsalpha}

\end{document}